\documentclass[11pt]{article}
\usepackage{graphicx}
\usepackage{tikz}
\usepackage{amssymb}
\usepackage{amsmath,amsthm}
\usepackage{hyperref}
\usepackage{verbatim}
\usepackage{amsfonts}
\usepackage[T1]{fontenc}
\usepackage{multicol}
\usepackage{color}
\hypersetup{
    colorlinks=true,
    citecolor=red,
    linkcolor=blue,
    filecolor=magenta,
    urlcolor=cyan,
    pdftitle={Combinatorial Invariants for Metacyclic Groups},
    pdfauthor={Mattia Brescia and Eshita Mazumdar},
}

\newtheorem{thm}{Theorem}[section]
\newtheorem{prop}{Proposition}[section]
\newtheorem{cor}[thm]{Corollary}

\newtheorem{lem}[thm]{Lemma}
\newtheorem{conj}[thm]{Conjecture}

\newtheorem{rem}[prop]{Remark}
\newtheorem{op}[thm]{Open Problem}

\newcommand{\no}{\noindent}
\newcommand{\bdot}{\mathbin{\boldsymbol{\cdot}}}

\title{Combinatorial Invariants for Metacyclic Groups}
\author{Mattia Brescia\footnote{Corresponding author. Department of Mathematics, University of Naples Federico II, Italy. Email:  mattia.brescia(at)unina.it}, and Eshita Mazumdar\footnote{Mathematical and Physical Sciences, Ahmedabad University, India. Email: eshita.mazumdar(at)ahduni.edu.in}}
\date{}

\begin{document}
\begingroup
\allowdisplaybreaks
\maketitle
\no
\begin{abstract}
We study product-one invariants of finite metacyclic groups, with particular attention to metacyclic $p$-groups. For semidirect products $C_n\rtimes C_m$ such that $\operatorname{ord}_n(\ell)=m$ and $\operatorname{ord}_q(\ell)=m$ for every prime divisor $q$ of $n$, we give a short proof of the known formula for the small Davenport constant and determine the exact value of the $\eta$-constant. When $n=p$ is prime, we use Vosper's theorem to obtain a growth result for products of non-trivial conjugacy classes contained in $C_p$, from which we derive a block-saturation principle for product-one sequences.

For split metacyclic $p$-groups, we determine the small Davenport constant in terms of the order and the exponent, and obtain bounds for the $\eta$-constant and the Erd\H{o}s--Ginzburg--Ziv constant. The Davenport bound is also extended to finite nilpotent split metacyclic groups. Our approach relies in part on a general quotient-lifting principle for the three invariants $\mathsf d$, $\eta$, and $\mathsf s$. For non-cyclic metacyclic $p$-groups, this principle applies through the Frattini quotient, which is isomorphic to $C_p^2$. We further establish upper and lower bounds for the Erd\H{o}s--Ginzburg--Ziv constant of arbitrary metacyclic $p$-groups, and show that both bounds are best possible within this class.

The lifting method also provides a short proof of a known theorem of Han and Zhang for groups $C_m\ltimes C_n$ with $m\mid n$. Finally, we determine the Olson constant of the dihedral group $D_{2p}$, for every prime $p>3$, in terms of the Olson constant of $C_p$.
\end{abstract}

\textbf{Keywords:} Product-one sequences, Small Davenport Constant, $\eta$-Constant, Erd\H{o}s--Ginzburg--Ziv Constant, Metacyclic $p$-groups, Olson Constant

2020 Mathematics Subject Classification: Primary 11B75; Secondary 20D60, 11P70.

\section{Introduction}
For given positive integers $m\le n$, we denote the set $\{m,m+1,\ldots,n\}$ by $[m,n]$. Let $(G,\cdot)$ be a finite group, written multiplicatively with identity element $1$, and let $(\mathcal{F}(G),\bdot)$ denote the free abelian monoid generated by $G$. By a sequence $S$ over $G$, we mean an element $S\in\mathcal{F}(G)$. Typically, for a sequence $S=g_1\bdot g_2\bdot\dotsc\bdot g_\ell$ over $G$, the non-negative integer $|S|=\ell$ denotes its length. For two sequences $S=g_1\bdot\dotsc\bdot g_{\ell_1}$ and $T=h_1\bdot\dotsc\bdot h_{\ell_2}$, we write $S\bdot T:=g_1\bdot\dotsc\bdot g_{\ell_1}\bdot h_1\bdot\dotsc\bdot h_{\ell_2}$. For any $g\in G$, the multiplicity of $g$ in $S$ is
$$
\mathsf{v}_g(S):=\bigl|\{i\in[1,\ell]:g_i=g\}\bigr|.
$$
A sequence $T$ over $G$ is a \textit{subsequence} of $S$, denoted by $T\mid S$, if $\mathsf{v}_g(T)\leq\mathsf{v}_g(S)$ for every $g\in G$. If $\mathsf{v}_g(S)\leq1$ for every $g\in G$, then $S$ is called squarefree. The notation $g_1^{[n_1]}\bdot\dotsc\bdot g_k^{[n_k]}$ means that $g_i$ occurs $n_i$ times. If $T\mid S$, we write $S\bdot T^{[-1]}$ for the sequence obtained from $S$ by removing the terms of $T$. For a non-empty subset $A\subseteq G$, we denote by $S(A)$ the subsequence of $S$ consisting of the terms that belong to $A$.

For a sequence $S=g_1\bdot\dotsc\bdot g_k$ over $G$, we define
$$
\pi(S):=\left\{\prod_{i=1}^k g_{\sigma(i)}:\sigma\in\mathfrak{S}_k\right\},\qquad
\Pi(S):=\bigcup_{\substack{T\mid S\\ |T|\ge1}}\pi(T),
$$
where $\mathfrak{S}_k$ is the symmetric group on $k$ letters. For $r\in[1,k]$, we also put
$$
\Pi_{\le r}(S):=\bigcup_{\substack{T\mid S\\ 1\le |T|\le r}}\pi(T),\qquad
\Pi_r(S):=\bigcup_{\substack{T\mid S\\ |T|=r}}\pi(T).
$$
We say that $S$ is a \textit{product-one sequence} if $1\in\pi(S)$, and that it is \textit{product-one free} if $1\notin\Pi(S)$. A {\it minimal product-one sequence} is a product-one sequence which cannot be partitioned into two non-trivial product-one subsequences.

The combinatorial invariants considered in this paper are the following:
\begin{itemize}
    \item the \textit{small Davenport constant}
    $$
    \mathsf{d}(G):=\max\bigl\{|S|:S\in\mathcal{F}(G)\text{ is product-one free}\bigr\};
    $$
    equivalently,
    $$
    \mathsf{d}(G)+1=\min\bigl\{k:|S|\ge k\text{ implies }1\in\Pi(S)\bigr\};
    $$
    \item the \textit{$\eta$-constant}
    $$
    \eta(G):=\min\bigl\{k:|S|\ge k\text{ implies }1\in\Pi_{\le\exp(G)}(S)\bigr\};
    $$
    \item the \textit{Erd\H{o}s--Ginzburg--Ziv constant}, or \textit{EGZ constant},
    $$
    \mathsf{s}(G):=\min\bigl\{k:|S|\ge k\text{ implies }1\in\Pi_{\exp(G)}(S)\bigr\};
    $$
    \item the \textit{Olson constant}
    $$
    \mathsf{Ol}(G):=\min\bigl\{k:\text{every squarefree sequence }S\text{ with }|S|\ge k\text{ satisfies }1\in\Pi(S)\bigr\}.
    $$
\end{itemize}

Recall that a finite group $G$ is \textit{metacyclic} if it has a cyclic normal subgroup $C$ such that $G/C$ is cyclic. A standard presentation is
\begin{equation}\label{eq1}
G=\langle x,y\mid y^n=1,\ x^m=y^t,\ y^x=y^\ell\rangle,
\end{equation}
where $\gcd(\ell,n)=1$ and the parameters satisfy the compatibility conditions
$$
\ell^m\equiv1\pmod n,\qquad t(\ell-1)\equiv0\pmod n.
$$ 
We refer to Hempel~\cite{Hemp} for the classification and the standard normalizations of these presentations. In the split case one may take $t=0$, and we shall usually write $G=B\ltimes C$, where $B$ and $C$ are cyclic, $C\trianglelefteq G$, and $B\cap C=\{1\}$. Our arguments for metacyclic $p$-groups use mainly this extension structure rather than the full parameter classification.

For the split groups appearing in connection with Bass' conjecture, we use the notation
$$
G_{m,n,\ell}=\langle x,y\mid x^m=y^n=1,\ y^x=y^\ell\rangle,
$$
where $\gcd(\ell,n)=1$ and $\ell^m\equiv1\pmod n$. In 2007, Bass~\cite{B} proposed, among other statements, the following conjectural formula for the small Davenport constant.

\begin{conj}\label{J.Bass}
Let $m\ge1$ and $n\ge2$ be integers and let $\ell$ be an integer such that $\operatorname{ord}_n(\ell)=m$. Consider the split metacyclic group
$$
G=\langle x,y\mid x^m=y^n=1,\ y^x=y^\ell\rangle.
$$
Then
$$
\mathsf{d}(G)=m+n-2.
$$
\end{conj}

In the same work, Bass determined the small Davenport constant for several basic non-abelian families, including dihedral and dicyclic groups, and proved the conjectural formula for groups of type $C_p\rtimes C_m$ with $p$ prime and $\operatorname{ord}_p(\ell)=m$. Brochero Mart\'inez and Ribas~\cite{BMR} subsequently classified the extremal product-one free sequences in this family. More recently, Wang and Zhao~\cite{WZ} proved Conjecture~\ref{J.Bass} under the additional assumption that $\operatorname{ord}_q(\ell)=m$ for every prime divisor $q$ of $n$, while Babu et al.~\cite{KBGS} obtained further cases of Bass' conjecture together with results on related Davenport-type invariants. In Section~3 we revisit the Wang--Zhao case: we give a short proof of
$$
\mathsf d(C_n\rtimes C_m)=n+m-2
$$
under the assumption that $\operatorname{ord}_q(\ell)=m$ for every prime divisor $q$ of $n$, and deduce
$$
\eta(C_n\rtimes C_m)=n+m-1.
$$
For $n=p$, we also prove the following product-growth bound. If $H=C_p$ and $\mathcal C_1,\ldots,\mathcal C_k$ are non-trivial conjugacy classes contained in $H$, not necessarily distinct, then
$$
|\mathcal C_1\cdots \mathcal C_k|\ge \min\{p,km\}.
$$
Vosper's theorem improves the bound obtained by repeated applications of Cauchy--Davenport. Since $m\mid p-1$, any $(p-1)/m+1$ such classes multiply to the whole of $H$. This yields a block-saturation statement for product-one sequences modulo $H$.

Han and Zhang~\cite{HZ} determined the small Davenport constant for groups $C_m\ltimes C_{mn}$ with arbitrary action. In Section~4 we give a short proof of the formula when $m\mid n$; the argument anticipates our use of quotients isomorphic to $C_p^2$.

Our main results concern metacyclic $p$-groups. We first determine the small Davenport constant of every split metacyclic $p$-group in terms of its order and exponent.

\begin{thm}\label{split}
Let $p$ be a prime and let $G$ be a non-trivial split metacyclic $p$-group. Then
$$
\mathsf{d}(G)=\exp(G)+\frac{|G|}{\exp(G)}-2.
$$
\end{thm}

For the metacyclic family $G_2(\alpha,\beta)$ considered by Godara, Joshi, and Mazumdar~\cite{GJM}, with $p$ odd, they proved that $\mathsf{d}(G)+1=p^\alpha+p^\beta-1$. Theorem~\ref{split} recovers this value and extends the formula to all split metacyclic $p$-groups, including $p=2$.

The lifting argument also gives the following nilpotent extension. Recall that a finite nilpotent split metacyclic group is the direct product of its Sylow subgroups, each of which is again split metacyclic; conversely, a direct product of split metacyclic Sylow subgroups is split metacyclic because the relevant cyclic factors have pairwise coprime orders.

\begin{thm}\label{Corsplit}
Let $G$ be a finite nilpotent split metacyclic group, and write
$$
G=H_1\times\cdots\times H_t
$$
for its Sylow decomposition, where $H_i$ is a split metacyclic $p_i$-group and the primes $p_1,\ldots,p_t$ are pairwise distinct. Then
$$
\mathsf{d}(G)\le\exp(G)+\frac{|G|}{\exp(G)}-2.
$$
Moreover, write $H_i=B_i\ltimes C_i$, where $B_i$ and $C_i$ are cyclic, $C_i\trianglelefteq H_i$, and $B_i\cap C_i=\{1\}$. Equality holds if either $\exp(H_i)=|C_i|$ for every $i\in[1,t]$, or $|H_i|/\exp(H_i)=|C_i|$ for every $i\in[1,t]$.
\end{thm}

For split metacyclic $p$-groups, we also obtain bounds for the $\eta$-constant and the EGZ constant.

\begin{thm}\label{split-three-invariants}
Let $p$ be a prime and let $G$ be a split metacyclic $p$-group. Put $E=\exp(G)$ and $Q=|G|/E$. Then
\begin{enumerate}
    \item $E+Q-1\le\eta(G)\le E+2Q-2$;
    \item $2E+Q-2\le\mathsf{s}(G)\le2E+2Q-3$.
\end{enumerate}
\end{thm}

In Section~2 we prove a quotient-lifting lemma for any normal subgroup $N\trianglelefteq G$, bounding the three invariants in terms of those of $N$ and $G/N$. For a non-cyclic metacyclic $p$-group, the Frattini quotient is isomorphic to $C_p^2$: if $G=\langle a,b\rangle$ with $\langle a\rangle\trianglelefteq G$, then
$$
\Phi(G)=\langle a^p,b^p\rangle,\qquad G/\Phi(G)\simeq C_p^2.
$$
We apply the lifting lemma through this quotient, which has the smallest possible order among non-cyclic quotients of $G$.

We also consider metacyclic $p$-groups that need not be split. At present, relatively little is known about $\mathsf{s}(G)$ for finite non-abelian groups. For cyclic groups one has $\mathsf{s}(G)=2|G|-1$ by the Erd\H{o}s--Ginzburg--Ziv theorem~\cite{EGZ}, while for rank-two abelian groups $C_{n_1}\times C_{n_2}$ with $n_1\mid n_2$ one has $\mathsf{s}(G)=2n_1+2n_2-3$; see~\cite{GH}. The non-abelian case has been investigated for special families such as dihedral and dicyclic groups~\cite{JZ}. For these groups, we prove the following bounds in terms of the order and exponent.

\begin{thm}\label{p-group}
Let $p$ be a prime and let $G$ be a finite metacyclic $p$-group. Put
$$
E=\exp(G),\qquad Q=\frac{|G|}{E}=p^q.
$$
Then
$$
2E+q(p-1)-1\leq\mathsf{s}(G)\leq2E+2Q-3.
$$
Both bounds are best possible over the class of finite metacyclic $p$-groups.
\end{thm}

Finally, we consider the Olson constant, the set-theoretic analogue of the small Davenport constant, which has been studied extensively for finite abelian groups. One of the earliest general results is due to Szemer\'edi~\cite{S}, who proved a conjecture of Erd\H{o}s and Heilbronn~\cite{EH}. Exact values are known for several important abelian families, including cyclic groups of prime order and various classes of $p$-groups and groups of small exponent; see, for instance,~\cite{Balandraud,OPS,Subocz}. Much less is known in the non-abelian case. For dihedral groups of order twice an odd prime, we prove the following formula.

\begin{thm}\label{Olson}
Let $p>3$ be a prime, let $G=D_{2p}$ be the dihedral group of order $2p$, and let $H=C_p$ be its cyclic subgroup of order $p$. Then
$$
\mathsf{Ol}(G)=\mathsf{Ol}(H)+1.
$$
\end{thm}

The restriction $p>3$ is necessary. Indeed, $\mathsf{Ol}(C_3)=2$, whereas the three reflections of $D_6$ form a product-one free subset: a product of one or three reflections is a reflection, while the product of two distinct reflections is a non-trivial rotation. Hence $\mathsf{Ol}(D_6)\ge4>3=\mathsf{Ol}(C_3)+1$.

The paper is organized as follows. Section~2 collects the auxiliary results, including the quotient-lifting lemma. Section~3 treats $C_n\rtimes C_m$ and the conjugacy-class growth and block-saturation results for $C_p\rtimes C_m$. Section~4 gives a short proof of the Han--Zhang theorem and our results on split metacyclic groups, including the nilpotent extension. Section~5 proves the EGZ bounds for arbitrary metacyclic $p$-groups, and Section~6 concerns the Olson constant of $D_{2p}$.

\section{Preliminaries}
This section collects the auxiliary results used in the proofs of our main theorems. For subsets $A,B$ of a group, we write $AB=\{ab:a\in A,\ b\in B\}$; when the ambient group is written additively, we use the corresponding notation $A+B$.

We shall use the following classical additive results in the case $n=p$. The first is the Cauchy--Davenport theorem, and the second is a theorem of Vosper~\cite{Vosper}.

\begin{thm}[Cauchy-Davenport~\cite{N}]\label{CD}
Let $p$ be a prime and let $A,B$ be non-empty subsets of $C_p$. Then
$$
|A+B|\ge \min\{p,|A|+|B|-1\}.
$$
\end{thm}

\begin{thm}[Vosper~\cite{Vosper}]\label{Vosper}
Let $p$ be a prime and let $A,B\subseteq C_p$ satisfy $|A|,|B|\ge2$. If
$$
|A+B|=|A|+|B|-1\le p-2,
$$
then $A$ and $B$ are arithmetic progressions with a common difference.
\end{thm}

We also use the following theorem of Scherk, originally proved as the solution of a problem of Moser~\cite{Scherk}.

\begin{lem}\label{MoserScherk}
Let $A$ and $B$ be finite subsets of an abelian group $G$ written additively, both containing the identity element $0$. If the only representation of $0=a+b$, with $a\in A$ and $b\in B$ implies that $a=b=0$, then
$$
|A+B|\ge \min\{|G|,|A|+|B|-1\}.
$$
\end{lem}

We next record a standard lower bound for the small Davenport constant under passage to a normal subgroup. We write the proof for the reader's convenience.

\begin{lem}\label{lower}
Let $G$ be a finite group and let $H$ be a normal subgroup of $G$. Then 
$$
\mathsf{d}(G)\ge \mathsf{d}(H)+\mathsf{d}(G/H).
$$
\end{lem}
\begin{proof}
Take product-one free sequences $S_1$ over $H$ and $T=(b_1H)\bdot\dotsc\bdot(b_rH)$ over $G/H$ of lengths $\mathsf{d}(H)$ and $r=\mathsf{d}(G/H)$, respectively, and choose representatives $b_i\in G$. Then $S=S_1\bdot b_1\bdot\dotsc\bdot b_r$ is product-one free. Thus $|S|=\mathsf{d}(H)+\mathsf{d}(G/H)$, proving the claim.
\end{proof}

We next recall a standard structural fact about the exponent of metacyclic $p$-groups. We include a proof to make explicit that the assertion holds for any pair of generators with one generating a normal cyclic subgroup.
\begin{lem}\label{metacyclic-exponent}
Let $p$ be a prime and let $H=\langle x,y\rangle$ be a finite metacyclic $p$-group with $\langle x\rangle\trianglelefteq H$. Then
$$
\exp(H)=\max\{|x|,|y|\}.
$$
\end{lem}
\begin{proof}
Write $|x|=p^u$, $|y|=p^v$, and $r=\max\{u,v\}$. Conjugation by $y$ induces an automorphism of $\langle x\rangle$ of order a power of $p$, so $x^y=x^\lambda$ for an integer $\lambda\equiv1\pmod p$ (if $x=1$, choose $\lambda=1$). Every element of $H$ can be written as $g=y^j x^i$. Putting $t=\lambda^j$, we obtain
$$
g^{p^r}=y^{jp^r}x^{i(1+t+\cdots+t^{p^r-1})}.
$$
The first factor is $1$. Moreover,
$$
1+t+\cdots+t^{p^r-1}
=\prod_{h=0}^{r-1}\bigl(1+t^{p^h}+\cdots+t^{(p-1)p^h}\bigr)
\equiv0\pmod{p^r},
$$
since each factor in the product is divisible by $p$. As $r\ge u$, the second factor is also $1$. Hence every element of $H$ has order dividing $p^r$, and the reverse inequality follows from the orders of $x$ and $y$.
\end{proof}

We will also need the following sharper lower bound for split metacyclic $p$-groups.

\begin{lem}\label{lower2}
Let $p$ be a prime and let $G=B\ltimes C$ be a split metacyclic $p$-group, where $B=\langle b\rangle$ and $C=\langle c\rangle\trianglelefteq G$ are cyclic and $B\cap C=\{1\}$. Then
$$
\mathsf{d}(G)\ge \exp(G)+\frac{|G|}{\exp(G)}-2.
$$
\end{lem}
\begin{proof}
By Lemma~\ref{metacyclic-exponent},
$$
\exp(G)=\max\{|B|,|C|\},
$$
and therefore
$$
\exp(G)+\frac{|G|}{\exp(G)}=|B|+|C|.
$$
Consider the sequence
$$
S=b^{[|B|-1]}\bdot c^{[|C|-1]}.
$$
If a non-trivial subsequence of $S$ were product-one, then, after passing to $G/C\simeq B$, the number of its terms equal to $b$ would be divisible by $|B|$. Since this number is strictly smaller than $|B|$, it must be zero. The subsequence would then consist only of copies of $c$, fewer than $|C|$, and could not be product-one. Thus $S$ is product-one free, and hence
$$
\mathsf{d}(G)\ge |S|=|B|+|C|-2=\exp(G)+\frac{|G|}{\exp(G)}-2.
$$
\end{proof}

We next collect the three classical values for $C_p^2$ that enter the lifting arguments.

\begin{lem}\label{Cp2constants}
Let $p$ be a prime. Then
$$
\mathsf{d}(C_p^2)=2p-2,\qquad
\eta(C_p^2)=3p-2,\qquad
\mathsf{s}(C_p^2)=4p-3.
$$
\end{lem}
The first equality is the classical formula for the small Davenport constant of finite abelian $p$-groups due to Olson~\cite{O1}; the second is also due to Olson~\cite{O2}; and the third is equivalent to Kemnitz's conjecture and follows from Reiher's theorem~\cite{ReiherKemnitz}.

The next lemma is the general lifting principle underlying the recursive arguments used later.

\begin{lem}\label{quotient-lift}
Let $G$ be a finite group, let $N\trianglelefteq G$, put $Q=G/N$, and let $e=\exp(Q)$. Then
\begin{enumerate}
    \item
    $$
    \mathsf{d}(G)+1\le e\,\mathsf{d}(N)+\max\{\mathsf{d}(Q)+1,\eta(Q)-e\}.
    $$
    Equivalently,
    $$
    \mathsf{d}(G)\le e\,\mathsf{d}(N)+\max\{\mathsf{d}(Q),\eta(Q)-e-1\}.
    $$
    \item If $\exp(G)=e\exp(N)$, then
    $$
    \eta(G)\le e(\eta(N)-1)+\eta(Q).
    $$
    \item If $\exp(G)=e\exp(N)$, then
    $$
    \mathsf{s}(G)\le e(\mathsf{s}(N)-1)+\mathsf{s}(Q).
    $$
\end{enumerate}
\end{lem}
\begin{proof}
For (1), put $r=\mathsf{d}(N)+1$ and
$$
L=e(r-1)+\max\{\mathsf{d}(Q)+1,\eta(Q)-e\}.
$$
Let $S$ be a sequence over $G$ of length at least $L$. As long as fewer than $r-1$ blocks have been removed, the image in $Q$ of the remaining sequence has length at least $e+\max\{\mathsf{d}(Q)+1,\eta(Q)-e\}\ge\eta(Q)$; hence one can extract a non-trivial subsequence of length at most $e$ whose product, in a suitable order, lies in $N$. After $r-1$ such blocks have been removed, at least $\mathsf{d}(Q)+1$ terms remain, so one can extract one further non-trivial subsequence whose product lies in $N$. Choosing an ordering in each block gives $r=\mathsf{d}(N)+1$ elements of $N$, some non-trivial subsequence of which has product $1$ in a suitable order. Replacing those elements by the corresponding blocks yields a product-one subsequence of $S$.

For (2), put $r=\eta(N)$ and let $S$ have length $e(r-1)+\eta(Q)$. Repeated application of the definition of $\eta(Q)$ to the image of the remaining sequence yields $r$ pairwise disjoint non-trivial subsequences $S_1,\ldots,S_r$, each of length at most $e$, whose terms can be ordered with products $x_1,\ldots,x_r\in N$. Since $r=\eta(N)$, some non-trivial subsequence of at most $\exp(N)$ of the $x_i$ has product $1$ in a suitable order. Replacing the selected $x_i$ by the corresponding blocks gives a product-one subsequence of $S$ of length at most $e\exp(N)=\exp(G)$.

For (3), put $r=\mathsf{s}(N)$ and let $S$ have length $e(r-1)+\mathsf{s}(Q)$. As in the preceding case, repeated application of the definition of $\mathsf{s}(Q)$ gives $r$ pairwise disjoint blocks $S_1,\ldots,S_r$, each of length exactly $e$, whose terms can be ordered with products $x_1,\ldots,x_r\in N$. Since $r=\mathsf{s}(N)$, exactly $\exp(N)$ of the $x_i$ have product $1$ in a suitable order. The corresponding blocks form a product-one subsequence of $S$ of length $e\exp(N)=\exp(G)$.
\end{proof}

Specializing the quotient-lifting lemma to the smallest non-cyclic elementary abelian quotient gives the form used throughout the rest of the paper.

\begin{cor}\label{Cp2lift}
Let $G$ be a finite group and let $N\trianglelefteq G$ satisfy $G/N\simeq C_p^2$. Then
$$
\mathsf{d}(G)\le p\,\mathsf{d}(N)+2p-2.
$$
If, in addition, $\exp(G)=p\exp(N)$, then
$$
\eta(G)\le p(\eta(N)-1)+3p-2
$$
and
$$
\mathsf{s}(G)\le p(\mathsf{s}(N)-1)+4p-3.
$$
\end{cor}
\begin{proof}
Apply Lemma~\ref{quotient-lift} and Lemma~\ref{Cp2constants}. For the first inequality, note that
$$
\max\{\mathsf{d}(C_p^2),\eta(C_p^2)-p-1\}
=\max\{2p-2,2p-3\}=2p-2.
$$
\end{proof}

For metacyclic $p$-groups, the normal subgroup occurring in Corollary~\ref{Cp2lift} is canonical.

\begin{lem}\label{metacyclic-frattini}
Let $p$ be a prime and let $G=\langle a,b\rangle$ be a finite non-cyclic metacyclic $p$-group with $C=\langle a\rangle\trianglelefteq G$ and $G/C=\langle bC\rangle$. Put
$$
N=\langle a^p,b^p\rangle.
$$
Then
$$
N=\Phi(G),\qquad G/N\simeq C_p^2,
$$
$N$ is metacyclic, and
$$
|N|=\frac{|G|}{p^2},\qquad \exp(N)=\frac{\exp(G)}p.
$$
If $G$ is split with $G=\langle b\rangle\ltimes\langle a\rangle$, then $N=\langle b^p\rangle\ltimes\langle a^p\rangle$ is split as well.
\end{lem}
\begin{proof}
Conjugation by $b$ induces an automorphism of the cyclic group $C$. Since the image of the $p$-group $\langle b\rangle$ in $\operatorname{Aut}(C/C^p)\simeq C_{p-1}$ is trivial, one has $[a,b]\in C^p$. As $G/C$ is cyclic, it follows that
$$
[G,G]\le C^p=\langle a^p\rangle.
$$
Consequently $N$ is normal: modulo $\langle a^p\rangle$ the generators $a$ and $b$ commute, so the subgroup generated by the image of $b^p$ is normal. Moreover, $G/N$ is generated by the commuting elements $aN$ and $bN$, each of order dividing $p$, and hence is elementary abelian of rank at most two.

Since $a^p,b^p\in\Phi(G)$, one has $N\le\Phi(G)$. If $G/N$ were cyclic, then the quotient $G/\Phi(G)$ would be cyclic as well, and the Burnside basis theorem would force $G$ to be cyclic, a contradiction. Hence $G/N\simeq C_p^2$. Since $G/N$ is elementary abelian, $\Phi(G)=G^p[G,G]\le N$, and therefore $N=\Phi(G)$.

The subgroup $C^p$ is cyclic and normal in $N$, while $N/C^p$ is generated by $b^pC^p$. The quotient calculation gives $|N|=|G|/p^2$. Finally, Lemma~\ref{metacyclic-exponent}, applied first to $G=\langle a,b\rangle$ and then to $N=\langle a^p,b^p\rangle$, gives
$$
\exp(N)=\max\{|a^p|,|b^p|\}=\frac1p\max\{|a|,|b|\}=\frac{\exp(G)}p.
$$
The last assertion is immediate in the split case because $\langle b^p\rangle\cap\langle a^p\rangle=\{1\}$.
\end{proof}

We recall the following known result.

\begin{lem}\label{universal-s-lower}
Let $G$ be a finite group and put $E=\exp(G)$. Then
$$
\mathsf{s}(G)\ge\eta(G)+E-1\ge \mathsf{d}(G)+E.
$$
\end{lem}
See, for instance, \cite[Section~1]{AMR}.

\begin{lem}\label{pgroup-D-lower}
Let $p$ be a prime and let $G$ be a finite $p$-group. Put $E=\exp(G)$, and write $|G|/E=p^q$. Then
$$
\mathsf{d}(G)\ge E+q(p-1)-1.
$$
\end{lem}
\begin{proof}
Since all element orders in a finite $p$-group are powers of $p$, there is an element $g\in G$ of order $E$. Put $H_0=\langle g\rangle$. Choose a maximal chain of subgroups
$$
H_0<H_1<\cdots<H_q=G.
$$
Since $|G:H_0|=p^q$, every inclusion in this chain has index $p$, and hence $H_{i-1}\trianglelefteq H_i$ and $H_i/H_{i-1}\simeq C_p$. Repeated application of Lemma~\ref{lower} gives
$$
\mathsf{d}(H_i)\ge \mathsf{d}(H_{i-1})+\mathsf{d}(C_p)=\mathsf{d}(H_{i-1})+p-1.
$$
As $\mathsf{d}(H_0)=E-1$, induction on $i$ yields $\mathsf{d}(G)\ge E+q(p-1)-1$.
\end{proof}

\section{The groups \texorpdfstring{$C_n\rtimes C_m$}{Cn semidirect Cm}}
We begin with a slightly more general result. Let $m,n\ge2$ and let $\ell$ be an integer such that
$$
\operatorname{ord}_n(\ell)=m
$$
and
$$
\operatorname{ord}_q(\ell)=m
$$
for every prime divisor $q$ of $n$. Consider
$$
G=G_{m,n,\ell}=\langle x,y\mid x^m=y^n=1,\ y^x=y^\ell\rangle,
$$
and put $H=\langle y\rangle\simeq C_n$. Then $G/H\simeq C_m$. The formula in the next theorem is contained in the more general result of Wang and Zhao~\cite{WZ}; its case $n=p$ was proved earlier by Bass~\cite{B}. We include a short proof because the rotation argument will also motivate the product-growth result below.

\begin{thm}\label{prime-kernel-d}
Let $m,n\ge2$, and let $\ell$ be an integer such that $\operatorname{ord}_n(\ell)=m$ and
$$
\operatorname{ord}_q(\ell)=m
$$
for every prime divisor $q$ of $n$. For
$$
G=\langle x,y\mid x^m=y^n=1,\ y^x=y^\ell\rangle,
$$
one has
$$
\mathsf d(G)=n+m-2.
$$
\end{thm}

\begin{proof}
The lower bound follows from Lemma~\ref{lower}, applied to $H\trianglelefteq G$:
$$
\mathsf d(G)\ge \mathsf d(H)+\mathsf d(G/H)=n+m-2.
$$
For the reverse inequality, let $S$ be a sequence over $G$ of length $n+m-1$, and suppose that $S$ is product-one free. Choose pairwise disjoint subsequences $T_1,\ldots,T_r$ of $S$ such that the product of the terms of each $T_i$ belongs to $H$, but no non-trivial proper subsequence of $T_i$ has product in $H$, and choose them maximally, so that the image in $G/H$ of the remaining sequence is product-one free. Since $\mathsf d(C_m)=m-1$,
$$
\sum_{i=1}^r |T_i|\ge |S|-(m-1)=n.
$$
For each $i$, every product of all terms of $T_i$ lies in $H$. We claim that
$$
|\pi(T_i)|\ge |T_i|.
$$
Fix an ordering $T_i=g_1\bdot\ldots\bdot g_t$ and put $h=g_1\cdots g_t\in H$. Since $S$ is product-one free, $h\ne1$.
Let $y^c\ne1$ be a non-trivial element of $H$ and choose a prime $q$ dividing $\operatorname{ord}(y^c)$. If $x^ay^b$ centralizes $y^c$, then $\ell^a\equiv1\pmod q$. Since $\operatorname{ord}_q(\ell)=m$, we get $m\mid a$, so $x^ay^b\in H$. Then the centralizer in $G$ of every non-trivial element of $H$ is exactly $H$ and in particular $C_G(h)=H$. Consider the $t$ cyclic rotations of $g_1\cdots g_t$. If two of their products were equal, a non-trivial proper consecutive block of the chosen ordering would centralize a non-trivial element of $H$, and hence its product would lie in $H$. Its image in $G/H$ would therefore be product-one, contradicting the minimality of the image of $T_i$. Thus the $t$ cyclic rotations have pairwise distinct products, proving the claim.

Put $A_i=\pi(T_i)\subseteq H$ and $A_i^*=A_i\cup\{1\}$. Since $S$ is product-one free, the identity has a unique representation
$$
1=a_1\cdots a_r,\qquad a_i\in A_i^*,
$$
namely $a_1=\cdots=a_r=1$: any other representation would give a product-one subsequence obtained by concatenating suitable orderings of the selected $T_i$. Iterating Lemma~\ref{MoserScherk}, we obtain
$$
|A_1^*\cdots A_r^*|\ge 1+\sum_{i=1}^r |A_i|
\ge1+\sum_{i=1}^r |T_i|\ge n+1,
$$
whereas $A_1^*\cdots A_r^*\subseteq H$ and $|H|=n$, a contradiction. Hence every sequence of length $n+m-1$ has a non-trivial product-one subsequence, and the upper bound follows.
\end{proof}

The same argument immediately determines the short product-one constant. Indeed, since $\operatorname{ord}_q(\ell)=m$ for every prime $q\mid n$, one has $m\mid q-1$ for every such $q$, and hence $\gcd(m,n)=1$. Therefore
$$
\exp(G)=mn.
$$

\begin{cor}\label{prime-kernel-eta}
Let $m,n\ge2$, and let $\ell$ be an integer such that $\operatorname{ord}_n(\ell)=m$ and
$$
\operatorname{ord}_q(\ell)=m
$$
for every prime divisor $q$ of $n$. For
$$
G=\langle x,y\mid x^m=y^n=1,\ y^x=y^\ell\rangle,
$$
one has
$$
\eta(G)=n+m-1.
$$
\end{cor}

\begin{proof}
The general inequality $\eta(G)\ge\mathsf d(G)+1$ and Theorem~\ref{prime-kernel-d} give the lower bound. Conversely, every sequence of length $n+m-1$ has a non-trivial product-one subsequence by Theorem~\ref{prime-kernel-d}, and its length is at most
$$
n+m-1\le mn=\exp(G).
$$
\end{proof}

We now return to the case $n=p$. Thus let
$$
G=G_{m,p,\ell}=\langle x,y\mid x^m=y^p=1,\ y^x=y^\ell\rangle,
$$
where $p$ is prime and $\operatorname{ord}_p(\ell)=m\ge2$, and put $H=\langle y\rangle\simeq C_p$. The extremal product-one free sequences in this family were classified by Brochero Mart\'inez and Ribas~\cite{BMR}.

We now exploit the full orbit structure of the non-trivial conjugacy classes contained in $H$. This is where Vosper's theorem yields a gain over repeated applications of Cauchy--Davenport.

\begin{thm}\label{prime-kernel-Vosper}
Let $p$ be a prime, let $m\ge3$, and let $\ell$ be an integer such that $\operatorname{ord}_p(\ell)=m$. Consider
$$
G=\langle x,y\mid x^m=y^p=1,\ y^x=y^\ell\rangle
\qquad\text{and}\qquad H=\langle y\rangle.
$$
Let $k\ge1$. If $\mathcal C_1,\ldots,\mathcal C_k$ are non-trivial conjugacy classes of $G$ contained in $H$, not necessarily distinct, then
$$
|\mathcal C_1\cdots \mathcal C_k|\ge \min\{p,km\}.
$$
In particular, if $q=(p-1)/m$, then
$$
\mathcal C_1\cdots \mathcal C_{q+1}=H.
$$
\end{thm}

\begin{proof}
Identify $H$ with the additive group $\mathbb F_p$. Put $L=\langle\ell\rangle\le\mathbb F_p^\times$, so $|L|=m$. Every non-trivial conjugacy class contained in $H$ has the form
$$
A_i=s_iL=\{s_i,s_i\ell,\ldots,s_i\ell^{m-1}\},\qquad s_i\in\mathbb F_p^\times.
$$
Thus it is enough to prove
$$
|A_1+\cdots+A_k|\ge\min\{p,km\}.
$$
We argue by induction on $k$. The case $k=1$ is immediate. Put
$$
A=A_1+\cdots+A_{k-1},\qquad B=A_k.
$$
If $|A|=p$, there is nothing to prove. Otherwise the induction hypothesis gives $|A|\ge(k-1)m$. Suppose, towards a contradiction, that
$$
|A+B|<\min\{p,km\}.
$$
By Cauchy--Davenport (Theorem \ref{CD}),
$$
|A+B|\ge \min\{p,|A|+m-1\}\ge \min\{p,km-1\}.
$$
If $km\ge p$, the preceding two inequalities force $km=p$, which is impossible because $m\mid p-1$. Hence $km<p$, and necessarily
$$
|A|=(k-1)m,\qquad |A+B|=|A|+|B|-1=km-1\le p-2.
$$
Vosper's Theorem \ref{Vosper} therefore implies that $B$ is an arithmetic progression.

We show that this is impossible. Since we are in the case $km<p$ and $k\ge2$, we have $m<p/2$. Suppose
$$
B=\{a,a+d,\ldots,a+(m-1)d\},\qquad d\ne0.
$$
Since $\ell B=B$, the same set is an arithmetic progression with common difference $\ell d$. For an arithmetic progression of length $3\le m<p/2$ in $C_p$, the common difference is determined up to sign: among the non-zero ordered differences of two elements, $d$ and $-d$ are precisely the two differences occurring $m-1$ times, while every other difference occurs fewer times. Hence $\ell d=\pm d$, so $\ell=\pm1$, contradicting $\operatorname{ord}_p(\ell)=m\ge3$. This proves the desired bound.

Finally, $(q+1)m=p-1+m\ge p$, so the first assertion gives $|\mathcal C_1\cdots \mathcal C_{q+1}|=p$, and hence $\mathcal C_1\cdots \mathcal C_{q+1}=H$.
\end{proof}

The preceding theorem has a direct interpretation for product-one blocks of maximal length in the quotient.

\begin{prop}\label{prime-kernel-blocks}
Let $p$ be a prime, let $m\ge3$, and let $\ell$ be an integer such that $\operatorname{ord}_p(\ell)=m$. Consider
$$
G=\langle x,y\mid x^m=y^p=1,\ y^x=y^\ell\rangle,
\qquad H=\langle y\rangle,
$$
and put $q=(p-1)/m$. Let $T_1,\ldots,T_{q+1}$ be sequences over $G$, each of length $m$, such that the image of every $T_i$ in $G/H\simeq C_m$ is a minimal product-one sequence. Then either one of the $T_i$ is product-one, or
$$
H\subseteq\pi(T_1\bdot\ldots\bdot T_{q+1}).
$$
In particular, in all cases $T_1\bdot\ldots\bdot T_{q+1}$ contains a non-trivial product-one subsequence; if none of the individual blocks is product-one, their full concatenation can be ordered to have product $1$.
\end{prop}

\begin{proof}
Assume that none of the $T_i$ is product-one. A minimal product-one sequence of length $m$ in $C_m$ consists of $m$ copies of a generator. Hence, after fixing an ordering of $T_i$, its product is a non-trivial element $h_i\in H$, and the products of the $m$ cyclic rotations run through the whole conjugacy class $\mathcal C_i=h_i^G\subseteq H$. Thus
$$
\mathcal C_i\subseteq\pi(T_i)
$$
for every $i$. Theorem~\ref{prime-kernel-Vosper} gives
$$
\mathcal C_1\cdots \mathcal C_{q+1}=H.
$$
Choosing independently an ordering of each $T_i$ that realizes the required element of $\mathcal C_i$ and concatenating these orderings shows that
$$
H\subseteq\pi(T_1\bdot\ldots\bdot T_{q+1}).
$$
In particular, $1$ belongs to this product set.
\end{proof}

\begin{rem}\label{prime-kernel-CD-gap}
The use of Vosper in Theorem~\ref{prime-kernel-Vosper} is essential for the stated saturation threshold. Repeated Cauchy--Davenport alone gives only
$$
|\mathcal C_1\cdots \mathcal C_k|\ge\min\{p,k(m-1)+1\},
$$
which does not in general force $H$ to be filled when $k=(p-1)/m+1$.
\end{rem}

\section{Product-one invariants of split metacyclic groups}

\subsection{The small Davenport constant}

We begin with a known result of Han and Zhang~\cite{HZ}, proved here in the case $m\mid n$. The proof introduces the quotient by $C_p^2$ that will be used for metacyclic $p$-groups.

\begin{thm}\label{Dmeta}
Let $m,n$ be positive integers with $m\mid n$, and let $\ell$ be an integer such that $\gcd(\ell,n)=1$ and $\ell^m\equiv1\pmod n$. For the split metacyclic group
$$
G=\langle x,y\mid x^m=y^n=1,\ y^x=y^\ell\rangle,
$$
one has
$$
\mathsf{d}(G)=m+n-2.
$$
\end{thm}

\begin{proof}
Let $C=\langle y\rangle\simeq C_n$. Since $C\trianglelefteq G$ and $G/C\simeq C_m$, Lemma~\ref{lower} gives
$$
\mathsf{d}(G)\ge \mathsf{d}(C_n)+\mathsf{d}(C_m)=n+m-2.
$$
It remains to prove the reverse inequality. We argue by induction on $|G|$. If $m=1$, then $G\simeq C_n$ and the assertion is immediate. Assume that $m>1$, let $p$ be the smallest prime divisor of $m$, and write
$$
m=ph,\qquad n=pk.
$$
Since $m\mid n$, we have $h\mid k$. Moreover, reducing the relation $\ell^m\equiv1\pmod n$ modulo $p$ gives $\ell^m\equiv1\pmod p$. Hence $\operatorname{ord}_p(\ell)$ divides both $m$ and $p-1$. By the minimality of $p$, one has $\gcd(m,p-1)=1$, and therefore
$$
\ell\equiv1\pmod p.
$$
Consider the homomorphism
$$
\varphi:G\longrightarrow C_p\times C_p,
$$
defined by sending $x$ and $y$ to the two standard generators. The defining relations of $G$ are respected because $p\mid m$, $p\mid n$, and $\ell\equiv1\pmod p$. Put
$$
K=\langle x^p,y^p\rangle.
$$
Since $\langle y^p\rangle$ is normalized by $x^p$, every element of $K$ has the form $x^{pa}y^{pb}$; thus $|K|=hk$. On the other hand, $K\subseteq\ker\varphi$ and $|\ker\varphi|=|G|/p^2=hk$. Consequently
$$
K=\ker\varphi\trianglelefteq G,\qquad G/K\simeq C_p^2.
$$
Furthermore, $K$ is again a split metacyclic group, generated by $x^p$ and $y^p$, of the form $C_h\ltimes C_k$. If $h=1$, then $K\simeq C_k$ and $\mathsf{d}(K)=k-1=h+k-2$; if $h>1$, the induction hypothesis applies and again gives
$$
\mathsf{d}(K)=h+k-2.
$$
Corollary~\ref{Cp2lift} now yields
$$
\mathsf{d}(G)\le p\,\mathsf{d}(K)+2p-2=p(h+k-2)+2p-2=m+n-2.
$$
Together with the lower bound, this proves the assertion.
\end{proof}

\begin{proof}[Proof of Theorem~\ref{split}]
Put $E=\exp(G)$ and $Q=|G|/E$. By Lemma~\ref{lower2}, $\mathsf{d}(G)\ge E+Q-2$. For the reverse inequality, we argue by induction on $|G|$. The cyclic case is immediate, and the first non-cyclic case $C_p^2$ follows from Lemma~\ref{Cp2constants}. Assume that $G$ is non-cyclic and not isomorphic to $C_p^2$. Write $G=C\rtimes B$, with $C$ and $B$ cyclic and $C\trianglelefteq G$, and put $N=\Phi(G)$. By Lemma~\ref{metacyclic-frattini},
$$
G/N\simeq C_p^2,\qquad \exp(N)=\frac Ep,\qquad |N|=\frac{|G|}{p^2},
$$
so that $|N|/\exp(N)=Q/p$. Moreover, $N$ is again split metacyclic. Hence, by induction,
$$
\mathsf{d}(N)=\frac Ep+\frac Qp-2.
$$
Corollary~\ref{Cp2lift} now gives
$$
\mathsf{d}(G)\le p\,\mathsf{d}(N)+2p-2=E+Q-2.
$$
Thus $\mathsf{d}(G)=E+Q-2$, as claimed.
\end{proof}

\begin{proof}[Proof of Theorem~\ref{Corsplit}]
Put
$$
E=\exp(G),\qquad Q=\frac{|G|}{E}.
$$
We first prove the upper bound by induction on $|G|$. If $G$ is cyclic, then $\mathsf{d}(G)=|G|-1=E+Q-2$. Suppose that $G$ is not cyclic. Since the Sylow factors have pairwise coprime orders, at least one $H_i$ is non-cyclic. Fix such an index $i$, put $p=p_i$, and set
$$
N=\Phi(H_i)\times\prod_{j\ne i}H_j.
$$
By Lemma~\ref{metacyclic-frattini}, $G/N\simeq C_p^2$, while $N$ is again a finite nilpotent split metacyclic group. Since the factor orders are pairwise coprime,
$$
\exp(N)=\frac Ep,\qquad |N|=\frac{|G|}{p^2},\qquad \frac{|N|}{\exp(N)}=\frac Qp.
$$
The induction hypothesis and Corollary~\ref{Cp2lift} yield
$$
\mathsf{d}(G)\le p\,\mathsf{d}(N)+2p-2\le E+Q-2.
$$

It remains to prove the asserted equality cases. Put $B=B_1\times\cdots\times B_t$ and $C=C_1\times\cdots\times C_t$. Because the orders belonging to distinct indices are pairwise coprime, both $B$ and $C$ are cyclic, $G=B\ltimes C$, $C\trianglelefteq G$, and $B\cap C=\{1\}$. If $b$ and $c$ generate $B$ and $C$, respectively, then
$$
b^{[|B|-1]}\bdot c^{[|C|-1]}
$$
is product-one free: projection to $G/C\simeq B$ first forces the number of selected $b$'s to be zero, and then the order of $c$ forces the number of selected $c$'s to be zero. Hence $\mathsf{d}(G)\ge|B|+|C|-2$.

If $\exp(H_i)=|C_i|$ for every $i$, then $E=|C|$ and $Q=|B|$; if $|H_i|/\exp(H_i)=|C_i|$ for every $i$, then $E=|B|$ and $Q=|C|$. In either case $\mathsf{d}(G)\ge E+Q-2$, so equality holds in the upper bound.
\end{proof}

\subsection{The short product-one and EGZ constants for split metacyclic p-groups}

\begin{proof}[Proof of Theorem~\ref{split-three-invariants}]
We prove the two upper bounds simultaneously by induction on $|G|$. If $G$ is cyclic, then $Q=1$ and the two asserted upper bounds are equalities. The group $C_p^2$ is the first non-cyclic case, covered by Lemma~\ref{Cp2constants}. Assume that $G$ is non-cyclic and not isomorphic to $C_p^2$. Put $N=\Phi(G)$. By Lemma~\ref{metacyclic-frattini},
$$
G/N\simeq C_p^2,\qquad \exp(N)=\frac Ep,\qquad \frac{|N|}{\exp(N)}=\frac Qp,
$$
and $N$ is again split metacyclic. The already proved Theorem~\ref{split} gives $\mathsf d(G)=E+Q-2$.

By induction,
$$
\eta(N)\le\frac Ep+\frac{2Q}{p}-2.
$$
Corollary~\ref{Cp2lift} yields
$$
\eta(G)\le p(\eta(N)-1)+3p-2\le E+2Q-2.
$$
Since $\mathsf{d}(G)+1\le\eta(G)$, we also have $E+Q-1\le\eta(G)$.

Similarly, by induction,
$$
\mathsf{s}(N)\le\frac{2E}{p}+\frac{2Q}{p}-3.
$$
Corollary~\ref{Cp2lift} gives
$$
\mathsf{s}(G)\le p(\mathsf{s}(N)-1)+4p-3\le2E+2Q-3.
$$
Finally, Lemma~\ref{universal-s-lower} and $\mathsf{d}(G)=E+Q-2$ imply
$$
\mathsf{s}(G)\ge \mathsf{d}(G)+E=2E+Q-2.
$$
This proves the two assertions.
\end{proof}

\begin{rem}\label{sharpness-three-invariants}
The upper bounds for $\eta(G)$ and $\mathsf{s}(G)$ in Theorem~\ref{split-three-invariants} are best possible over the full class of split metacyclic $p$-groups; in particular, the upper bound for $\mathsf{s}(G)$ in Theorem~\ref{p-group} is best possible for the class of all metacyclic $p$-groups. Indeed, for the trivial action one has $G\simeq C_E\times C_Q$ with $Q\mid E$, and the classical rank-two formulas give
$$
\eta(G)=E+2Q-2,\qquad \mathsf{s}(G)=2E+2Q-3.
$$
Thus the lower bound $\mathsf{s}(G)\ge2E+Q-2$ cannot be an equality uniformly throughout the split metacyclic class.
\end{rem}

\begin{rem}\label{qwl-index-p}
Qu, Wang, and Li~\cite{QWL} prove that if a finite non-abelian group has a cyclic subgroup of index $p$, where $p$ is the smallest prime divisor of its order, then
$$
\eta(G)=\mathsf{d}(G)+1,\qquad \mathsf{s}(G)=\eta(G)+\exp(G)-1.
$$
Consequently, if $G$ is a non-abelian split metacyclic $p$-group with $Q=|G|/E=p$, then
$$
\eta(G)=E+p-1,\qquad \mathsf{s}(G)=2E+p-2=2E+Q-2.
$$
Hence the lower bounds in Theorem~\ref{split-three-invariants} are attained in this subclass. Since here $Q=p$, the lower bound in Theorem~\ref{p-group} is attained as well.
\end{rem}

\begin{op}\label{op-split-eta-s}
Determine $\eta(G)$ and $\mathsf{s}(G)$ for non-abelian split metacyclic $p$-groups beyond the index-$p$ case of Remark~\ref{qwl-index-p}. In particular, characterize the groups for which
$$
\eta(G)=\mathsf{d}(G)+1=E+Q-1
$$
and those for which
$$
\mathsf{s}(G)=\eta(G)+E-1.
$$
In the simultaneous extremal case these identities give $\mathsf{s}(G)=\mathsf{d}(G)+E=2E+Q-2$.
\end{op}

\begin{rem}
Theorem~\ref{split-three-invariants} shows that determining $\eta(G)$ is an essential intermediate problem. The possible range
$$
E+Q-1\le\eta(G)\le E+2Q-2
$$
has width $Q-1$, and the corresponding range for $\mathsf{s}(G)$ has the same width.
\end{rem}

\section{EGZ bounds for metacyclic p-groups}

\begin{proof}[Proof of Theorem~\ref{p-group}]
Put $E=\exp(G)$ and $Q=|G|/E=p^q$. Lemmas~\ref{universal-s-lower} and~\ref{pgroup-D-lower} give
$$
\mathsf{s}(G)\ge \mathsf{d}(G)+E\ge2E+q(p-1)-1,
$$
which proves the lower bound.

We prove the upper bound by induction on $|G|$. If $G$ is cyclic, then $Q=1$ and
$$
\mathsf{s}(G)=2E-1=2E+2Q-3.
$$
Assume that $G$ is non-cyclic and put $N=\Phi(G)$. By Lemma~\ref{metacyclic-frattini}, $N$ is again metacyclic and
$$
G/N\simeq C_p^2,\qquad |N|=\frac{|G|}{p^2},\qquad \exp(N)=\frac Ep.
$$
Consequently
$$
\frac{|N|}{\exp(N)}=\frac Qp.
$$
If $N$ is trivial, then $G\simeq C_p^2$ and the desired upper bound is the classical equality $\mathsf{s}(C_p^2)=4p-3$. Otherwise the induction hypothesis gives
$$
\mathsf{s}(N)\le\frac{2E}{p}+\frac{2Q}{p}-3.
$$
Applying Corollary~\ref{Cp2lift}, we obtain
$$
\mathsf{s}(G)
\le p(\mathsf{s}(N)-1)+4p-3
\le p\left(\frac{2E}{p}+\frac{2Q}{p}-4\right)+4p-3
=2E+2Q-3.
$$
This proves the theorem.
\end{proof}

\begin{rem}\label{p-group-sharpness}
Both estimates in Theorem~\ref{p-group} are attained on infinite families. For the upper bound, let $E$ and $Q$ be powers of $p$ with $Q\mid E$ and take
$$
G=C_E\times C_Q.
$$
Then $\exp(G)=E$, $|G|/E=Q$, and the rank-two formula gives
$$
\mathsf{s}(G)=2E+2Q-3.
$$
Thus the upper bound is pointwise best possible for every such pair $(E,Q)$. For the lower bound, it is already sharp on a non-abelian family with $Q=p$. For example, if $p$ is odd and $n\ge2$, let
$$
G_{p,p^n}=C_{p^n}\rtimes C_p,
$$
where a generator of $C_p$ acts on $C_{p^n}$ by the automorphism $a\mapsto a^{1+p^{n-1}}$. Then $G_{p,p^n}$ is a non-abelian split metacyclic $p$-group with $E=p^n$ and $Q=p$. By Remark~\ref{qwl-index-p},
$$
\mathsf{s}(G_{p,p^n})=2E+p-2=2E+(p-1)\log_pQ-1.
$$
Hence the lower bound is best possible as well.
\end{rem}

\section{The Olson constant for dihedral groups}

\begin{proof}[Proof of Theorem~\ref{Olson}]
Note that $G=G_{2p}=\langle x,y\mid x^2=y^p=1,\ y^x=y^{-1}\rangle$ and $H=\langle y\rangle$. Put $k=\mathsf{Ol}(H)$. If $A=\{y_1,\ldots,y_{k-1}\}\subseteq H$ is product-one free, then $A\cup\{x\}$ is product-one free in $G$: indeed, every product involving $x$ lies in $G\setminus H$, because a subsequence of $A\cup\{x\}$ contains $x$ at most once. Hence $\mathsf{Ol}(G)\ge k+1$.

We prove the reverse inequality. Let $S$ be a square-free sequence over $G$ with $|S|=k+1$, and suppose, towards a contradiction, that $S$ has no non-trivial product-one subsequence. Write $S=S(H)S(N)$, where $N=G\setminus H$, and put $|S(H)|=u$ and $|S(N)|=v$. If $u\ge k$, then the definition of $k$ gives a product-one subsequence contained in $S(H)$, contrary to our assumption. Thus $u\le k-1$ and $v\ge2$.

Write the distinct elements of $S(N)$ as $xy^r$, with $r$ in a subset $R\subseteq \mathbb Z_p$ of cardinality $v$. If $v=p$, then $R=\mathbb Z_p$ and the four distinct reflections $x,xy,xy^3,xy^2$ have product $1$, since $(x\,xy)(xy^3\,xy^2)=y\,y^{-1}=1$. Hence $v<p$.

We claim that the elements of $S(N)$ can be ordered as $c_1,\ldots,c_v$ in such a way that no three consecutive exponents form an arithmetic progression in $\mathbb Z_p$. To see this, let $M$ be the number of ordered triples $(a,b,c)$ of distinct elements of $R$ satisfying $b-a=c-b$. For each ordered pair $(a,b)$ with $a\ne b$ there is at most one possible value $c=2b-a$, so $M\le v(v-1)$. In fact $M<v(v-1)$: equality would imply that $2b-a\in R$ for every distinct $a,b\in R$; fixing such a pair and iterating this relation would then give $a+n(b-a)\in R$ for every $n\in\mathbb Z_p$, hence $R=\mathbb Z_p$, contrary to $v<p$. For each fixed ordered arithmetic-progression triple $(a,b,c)$, there are exactly $(v-2)!$ orderings of $R$ in which $a,b,c$ occur consecutively, in this order: there are $v-2$ possible positions for the first term $a$, and the remaining $v-3$ elements can be arranged arbitrarily. Hence the number of orderings containing at least one consecutive arithmetic-progression triple is at most
$$
M(v-2)!<v(v-1)(v-2)!=v!.
$$
Since $R$ has exactly $v!$ orderings, there exists an ordering containing no three consecutive elements in arithmetic progression. Let $c_i=xy^{r_i}$ be the corresponding ordering of $S(N)$.

Set $h_i=c_ic_{i+1}\in H$ for $1\le i\le v-1$. First, $h_i\notin S(H)$ for every $i$: otherwise $c_{i+1}c_i h_i=1$ would be a product-one subsequence of $S$. Next, we claim that the elements $h_1,\ldots,h_{v-1}$ are pairwise distinct. Indeed, if $h_i=h_{i+1}$, then, writing $c_t=xy^{r_t}$, we obtain $r_{i+1}-r_i=r_{i+2}-r_{i+1}$, contrary to the choice of the ordering. If instead $h_i=h_j$ with $j\ge i+2$, then the four reflections involved are distinct and
$$
c_jc_ic_{i+1}c_{j+1}=1,
$$
again a contradiction. Thus
$$
T=S(H)\cdot h_1\cdots h_{v-1}
$$
is a square-free sequence over $H$ of length $u+v-1=k$. By the definition of $k$, it has a non-trivial product-one subsequence.

If this product-one subsequence uses only terms of $S(H)$, we obtain an immediate contradiction. Otherwise, let $I\subseteq\{1,\ldots,v-1\}$ be the set of indices of the selected terms $h_i$. Decompose $I$ into its maximal intervals of consecutive integers,
$$
I=[a_1,b_1]\,\dot\cup\cdots\dot\cup[a_t,b_t].
$$
Since $H$ is abelian and $c_i c_{i+1}=y^{r_{i+1}-r_i}$, the product over each interval simplifies by cancellation:
$$
\prod_{i=a_j}^{b_j}h_i=y^{r_{b_j+1}-r_{a_j}}=c_{a_j}c_{b_j+1}.
$$
The maximality of the intervals implies that the $2t$ reflections $c_{a_j},c_{b_j+1}$ are pairwise distinct. Hence, replacing in the product-one subsequence of $T$ each block $\prod_{i=a_j}^{b_j}h_i$ by the two original terms $c_{a_j},c_{b_j+1}$, and ordering first the selected terms from $S(H)$ and then these reflection pairs, gives a non-trivial product-one subsequence of the original sequence $S$. This contradiction proves $\mathsf{Ol}(G)\le k+1$, and therefore $\mathsf{Ol}(G)=\mathsf{Ol}(H)+1$.
\end{proof}

\section{Acknowledgements}
The first author gratefully acknowledges the support and hospitality of Prof. Mazumdar and of the School of Arts and Sciences during the preparation of this paper at Ahmedabad University, as well as the funding by the ``Programma di Scambi Internazionali con Università ed Istituti di Ricerca Stranieri per la Mobilità di Breve Durata di Docenti, Ricercatori e Studiosi''. The second author gratefully acknowledges the support of the Science and Engineering Research Board (SERB), India, under the ANRF Core Research Grant (File No.\ CRG/  2023/002698).

\endgroup
\end{document}